\documentclass[conference]{IEEEtran}
\usepackage{cite}
\usepackage{amsmath,amssymb,amsfonts}
\usepackage{algorithmic}
\usepackage{graphicx}
\usepackage{textcomp}
\usepackage{xcolor}
\usepackage{subfig}
\def\BibTeX{{\rm B\kern-.05em{\sc i\kern-.025em b}\kern-.08em
    T\kern-.1667em\lower.7ex\hbox{E}\kern-.125emX}}
\begin{document}

\title{Composed Physics- and Data-driven System Identification for Non-autonomous Systems in Control Engineering}

\author{\IEEEauthorblockN{Ricarda-Samantha Götte}
\IEEEauthorblockA{\textit{Heinz Nixdorf Institute} \\
\textit{Paderborn University}\\
Paderborn, Germany \\
rgoette@hni.upb.de}
\and
\IEEEauthorblockN{Julia Timmermann}
\IEEEauthorblockA{\textit{Heinz Nixdorf Institute} \\
\textit{Paderborn University}\\
Paderborn, Germany \\
julia.timmermann@hni.upb.de}

}

\maketitle

\begin{abstract}
	In control design most control strategies are model-based and require accurate models to be applied successfully. Due to simplifications and the model-reality-gap physics-derived models frequently exhibit deviations from real-world-systems. Likewise, purely data-driven methods often do not generalise well enough and may violate physical laws. Recently Physics-Guided Neural Networks (PGNN) and physics-inspired loss functions separately have shown promising results to conquer these drawbacks. In this contribution we extend existing methods towards the identification of non-autonomous systems and propose a combined approach PGNN-L, which uses a PGNN and a physics-inspired loss term (-L) to successfully identify the system's dynamics, while maintaining the consistency with physical laws. The proposed method is demonstrated on two real-world nonlinear systems and outperforms existing techniques regarding complexity and reliability.
\end{abstract}


\section{Introduction}
\begin{figure}[!htb]
	\centering
	\subfloat[Residual Modelling]{\includegraphics[width=0.5\textwidth]{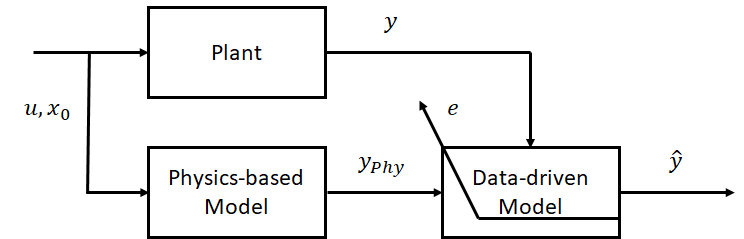}\label{fig:ResidualModeling}}
	\hfill
	\subfloat[Extended Residual Modelling]{\includegraphics[width=0.5\textwidth]{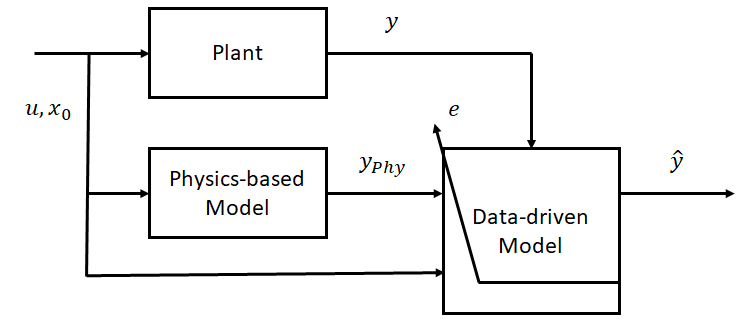}\label{fig:ExtendedModeling}}
	\caption{Schemes of approaches for combining data-driven and physics-based models: (\ref{fig:ResidualModeling}) The data-driven model approximates the deviation between the physics-based model's predictions $y_{phy}$ and the system's measurements $y$ without diving deeper into the error dynamics. (\ref{fig:ExtendedModeling}) The data-driven model takes the excitation $u$ and the initial state $x_0$ additionally into account, which enables it to recognise the deviations of the physics-based model and utilise these to model the system's dynamics itself.}
\end{figure}

In many applied engineering fields, such as control engineering, identifying suitable models is a crucial step to predict a system's behaviour and be able to control it with a model-based method. Deriving these models by classical scientific approaches, e.g. by formulating the governing equations, preserves the physical principles and therefore enables a physical understanding of the considered system. Nonetheless, these models often exhibit significant deviations from observed quantities due to the model-reality-gap or unmodelled partial dynamics, which cannot be decreased by simple parameter identification techniques.\newline
In recent years this challenge is addressed by taking measurement data into account during the modelling process. This allows a deeper insight into the underlying dynamics \cite{Brunton.2016}. However, data may be limited due to costly experiments, e.g. in automotive test benches, or may be incomplete. Additionally, purely data-driven models lack of physical context. Hence, they may not be capable of describing the system's dynamics beyond the scope of the used training data. Eventually, this may cause faulty predictions which can harm expensive hardware or increase risks of injury. Thus, currently a field of research on hybrid modelling emerges which seeks to combine physics-based and data-driven approaches towards a more sophisticated modelling \cite{Karpatne.2017}.\newline
From our point of view there are three perspectives how to use machine learning approaches and physics-based methods or vice versa in conjunction. When computing power had been quite low, physics-based modelling with data-driven aid has been proposed, for instance by extending a physical model with a neural network for parameter estimation \cite{Psichogios.1995}. With advancements in storage capacities and computation efficiency, the amount of accessible and processible data increased so that data-driven techniques became popular and physical models rather were utilised for pre-training \cite{Willard.2020}. Due to the awareness of potential downsides by considering solely one perspective over the other and the acknowledgement of the comprehensive physical knowledge a third perspective arises. Its aim is to assemble physics-based methods and machine learning approaches in a systematic, synergetic way to exploit the domain-specific understanding we usually have and the information we observe due to data \cite{Karpatne.2017, Willard.2020}. This approach emerges among different scientific disciplines, e.g. in hydrology, climatology, geology, fluid dynamics, and control engineering \cite{Karpatne.2017b, Willard.2020, Antonelo.2021}, which illustrates the great promise this field seems to hold for obtaining a deeper understanding of scientific phenomena.\newline        
Motivated by the model-reality-gap one of the most popular, yet simple hybrid approaches is \textit{Residual Modelling} \cite{Willard.2020}, which generally strives to approximate the error between physically derived model's outputs and system's measurement data by a data-driven method. Common techniques of modelling the error are least-squares optimisation or neural networks. The advantage of this approach is the simplicity, still it does not enable a deeper insight into possibly unmodelled dynamics, which occur as the deviation between model and measurement, or hold any physical meaning as it just models the deviation from the model to the real-world system. A visualisation of this basic idea is shown in figure \ref{fig:ResidualModeling}.
Extending the notion of residual modelling leads to an approach which enhances the simple serial connection of model and data-driven error model and rather takes the error itself directly into account. Recently, many data-driven techniques rely not only on measurement data but also on simulated data from an existing model. This additional data from simulation is explicitly taken into account depending on the structure of the method. As result, the model does not only approximate the error but also learns the relation between inputs, simulated outputs and measurements as it is shown in figure \ref{fig:ExtendedModeling}. A very popular technique are Physics-Guided or Physics-Informed Neural Networks (PGNNs or PINNs respectively), which consider the simulated model's outputs by augmenting the input vector or by using it as loss function within the learning procedure \cite{Karpatne.2017,Pathak.2018,Willard.2020, Yu.2020,Yang.2020}. Assuming that the prior physical model is already to somewhat extent accurate, in the first approach the neural network is guided by the additional input and mediates when deviations occur. This allows much lower training cost, generalisability by using the simulation model, which is applicable to various initial values, and a higher likelihood of physical consistency. However, this comes with the drawback of the partial loss of physical interpretability. To mitigate this, some research approaches include domain specific knowledge inside the learning process, e.g. by incorporating bounds or specific physical relations, which should guarantee the network's consistency with physical laws \cite{Muralidhar.2018, Karpatne.2017,  Willard.2020, Cranmer.2020, Greydanus.2019, Raymond.2021, Rueden.2019}.\newline		
In \cite{Karpatne.2017} the authors utilise a PGNN for the estimation of lake temperatures depending on the depth of a lake based on a general lake model. Additionally, they propose a physically motivated term added to the traditional loss function, which ensures physical consistency by respecting the relations of temperature and density as well as density and depth. This allows the learning procedure to respect physical relations and provides physically meaningful solutions. In the authors' work it can be recognised that this method provides much lower training error and ensures a higher level of physical correctness compared to other data-driven methods. Similar to the notion of incorporating physical relations as an additional loss term, there are also approaches that exclusively rely on a physical loss function as error function evaluated on estimated and targeted outputs, such as the Hamiltonian function \cite{Greydanus.2019}, the Lagrangian formulation of a system \cite{Cranmer.2020, Lutter.2015} or the governing equations \cite{Antonelo.2021}. If the prior knowledge is accurate enough, these methods propose good predictions. However, e.g. in case of the Hamiltonian method, the approach is restricted to conservative systems which are not likely in practical applications. Moreover, one needs explicit knowledge about the structure and governing terms of the considered system to formulate the Lagrangian or the ODEs, which is not realistic from a practical perspective.    
\newline
Other hybrid approaches include imposing physics into the architecture of neural networks as intermediate states \cite{Daw.2019}, but as stated before these methods require a good knowledge of the considered system's structure to be applied successfully. Moreover, several techniques focus on neural ordinary equations or even extend that notion towards a neural operator \cite{Chen.2018,Li.,Li.2021}. Depending on the considered system, sometimes a combination of multiple data-driven techniques might be useful, e.g. the formulation of Bayesian neural network as PGNN to solve partial differential equations \cite{Yang.2020}.\newline
In this contribution we extend the framework of PGNNs towards the system identification of non-autonomous systems and ensure physically consistent predictions by establishing a physically motivated loss term within the learning procedure. The proposed approach PGNN-L consists of a PGNN, whose input is augmented with the control input and the simulated output of a physical model, which describes the a priori presumed dynamics and guides the learning process as external input. This is based on the assumption that the prior model covers qualitatively the majority of the considered system's dynamics, but is not completely reliable in terms of unidentified parameters and partial dynamics. In contrast to \cite{Antonelo.2021}, which utilises the presumed dynamics within the learning procedure, we impose similar to other researchers \cite{Karpatne.2017, Muralidhar.2018, Raymond.2021} loss functions regarding qualitative information, such as the energy conservation principle, and investigate the weighting of multiple loss terms which has been examined very recently for the first time \cite{Rohrhofer.2021}. The proposed method guarantees a higher level of physical consistency and outperforms purely physics- or data-driven methods.\newline
The paper is organised as follows. Section \ref{sec:modelling} summarises classical mathematical or physical modelling approaches and gives an overview of data-driven methods for identifying mechatronic systems, whereupon our contribution is presented in detail. In section \ref{sec:experiments} the hybrid modelling approach is demonstrated and evaluated on two nonlinear systems. The paper concludes with a short discussion of the results in section \ref{sec:conclusion}.

\section{Modelling approaches}\label{sec:modelling}

A continuous, nonlinear, non-autonomous system is considered
\begin{align}\label{eq:system}
\begin{split}
\dot{x}&=f(x,u,p),\\
y&=g(x,u,p),
\end{split}
\end{align}
where $x\in\mathbb{R}^n$ describes the state, $u\in\mathbb{R}^m$ the control input, $y\in\mathbb{R}^l$ the measured output, $p\in\mathbb{R}^q$ the parameters of the system and $f:\,\mathbb{R}^n\times\mathbb{R}^m\times\mathbb{R}^q\mapsto\mathbb{R}^n$ and $g:\nobreak\,\mathbb{R}^n\times\mathbb{R}^m\times\mathbb{R}^q\mapsto\nobreak\mathbb{R}^l$ are nonlinear functions, describing the system's dynamics. In control engineering, we wish to model the plant \eqref{eq:system} and identify the parameters $p$ to proceed with the model-based control design, including the synthesis, analysis and simulation of each design step, e.g. the establishment of a controller. In this work we focus on the aforementioned first step of control design, the system and parameter identification.

\subsection{Mathematical Modelling}\label{subsec:mathematical}
Traditionally, real-world dynamics are characterised by a set of governing equations, derived from physical laws. Depending on the considered system, either Newton's axioms or Euler-Lagrange formalism are utilised for formulating continuous, potentially nonlinear models which describe the system's dynamics by ordinary differential equations (ODEs).\newline
In this work we assume that the considered system \eqref{eq:system} can be captured by the following system of ODEs
\begin{align}\label{eq:physModel}
\begin{split}
\dot{x}_{phy}&=f_{phy}(x_{phy},u,\tilde{p}),\\
y_{phy}&=g_{phy}(x_{phy},u,\tilde{p}),
\end{split}
\end{align}
in a qualitative manner which means it is not as accurate as we wish due to possibly unidentified parameters $\tilde{p}\in\mathbb{R}^q$ or unmodelled partial dynamics. This may be the case e.g. during the commissioning of a plant or when complex dynamical elements are neglected or unknown.\newline
In many fields, not only in control engineering, this model \eqref{eq:physModel} often is too complex in terms of dimension or nonlinearity to be utilised for any method, so that it needs to be simplified by linearisation techniques, such as Taylor expansion within concrete operating points. These models then can be applied in model-based methods, which require a fast and simple evaluation of the model to be able to predict, for instance model-predictive control, gain scheduling approaches and similar methods. The disadvantage of those derived models is primarily the expense to deduce and store a model for each operating point, but more important the potential of errors due to the simplified model. Therefore, it is always desirable to obtain a good trade-off between a model's accuracy and its complexity.

\subsection{Data-driven Modelling}\label{subsec:data-driven}
Lately, the utilisation of data in science and engineering becomes increasingly popular. 
In modelling complex processes, data-driven techniques may offer the opportunity to optimise the process and save resources like time, energy or manpower, for instance throughout the set-up and deployment of a novel plant. During regular operation but also even before, when installing and testing a novel plant by applying small excitations, data $D\subset\mathbb{R}^{(l+m)\times N}$ can be collected by $N\in\mathbb{N}$ measurement samples of the plant's excitation and output, $u_{M,k}$ and $y_{M,k}$ respectively, and structured time $t(k)=t_k$ dependently in the following way
\begin{align}
\begin{split}
D &= \Big\{  \left( y (t_k), u (t_k)\right) \Big\} \\&= \Big\{ \left( y_M, u_M\right)_k \Big\}  \quad k = 0,1,\dots,N.
\end{split}
\end{align} 
Thereby, the samples need not necessarily belong to one trajectory, rather should cover various states and excitations in the phase and input space except if the plant is in the process of deployment. Using these samples, different methods can be applied to derive a data-driven model which approximates the system \eqref{eq:system} to the best possible accuracy. \newline
\textit{Least-squares optimisation} is a very popular approach among machine learning techniques, which attempts to fit measurement data to a certain model, e.g. a linear regression model, by reducing the least-squares error \cite{Bishop.2006}.\newline		
A similar approach, utilising least-squares optimisation, is called \textit{Sparse Identification of Nonlinear Dynamics with control (SINDYc)} and based on the assumption that many systems' dynamics can be described by only a small number of governing terms \cite{Brunton.2016}. Organising the measurement data $D$, derived from the system \eqref{eq:system}, into snapshot matrices
\begin{align}
\begin{split}
Y &= \begin{pmatrix}
\vline & \vline &  & \vline\\
y_{M,1} & y_{M,2} & \dots & y_{M,N-1}\\
\vline & \vline &\  & \vline
\end{pmatrix} \in\mathbb{R}^{l\times (N-1)},\\
Y' & =\begin{pmatrix}
\vline & \vline &  & \vline\\
y_{M,2} & y_{M,3} & \dots & y_{M,N}\\
\vline & \vline &  & \vline
\end{pmatrix} \in\mathbb{R}^{l\times (N-1)},\\
U &= \begin{pmatrix}
\vline & \vline &  & \vline\\
u_{M,1} & u_{M,2} & \dots & u_{M,N-1}\\
\vline & \vline & & \vline
\end{pmatrix} \in\mathbb{R}^{m\times (N-1)},
\end{split}
\end{align}
the aim is to find a best-fit approximation model $Y' = \Xi\Psi(Y,U)^T$ by storing and evaluating $\kappa$ initial guess terms in a library vector $\Psi(Y,U)\in\mathbb{R}^{(N-1)\times \kappa}$, which might describe the considered system's dynamics \eqref{eq:system} well. The method then seeks a sparse vector $\Xi\in\mathbb{R}^{l\times \kappa}$ such that the following holds for each row $i$ with $0<\lambda\ll 1$ as a model complexity parameter
\begin{equation}
\xi_i^* = \text{argmin}_{\xi_i} \frac{1}{2}\left| \left|  Y'_i-\xi_i\Psi(Y,U)^T\right| \right| _2^2 +\lambda \left|\xi_i\right|.
\end{equation}
This technique is usually considered as data-driven method as it can be utilised without any prior knowledge. However, if prior knowledge is available, SINDYc moves towards the perspective of hybrid modelling approaches by incorporating preassumed terms into the library vector. Advantages of SINDYc include the computational efficiency, its real-time capability, once $\Xi$ is determined, and the preservation of the system's physical interpretability by $\Xi$, indicating which terms are relevant for the description of system \eqref{eq:system}. Nonetheless, it relies on the assumption of a few governing terms, which is a huge drawback because it might not be true for many complex systems. Furthermore, its model is not capable of imposing any physical or operational limitations, restricting its applicability in many mechatronic systems.\newline
\textit{Neural Networks (NNs)} have become very successful during the last years in a variety of applications, e.g. learning control policies \cite{Chang.} or addressing the problem of non-measurable states \cite{Chen.2018b}. Data are fed into input layers, processed through a number of hidden layers with (nonlinear) activation functions and transformed to estimated outputs $\hat{Y}$, which are compared to labelled data $Y$ in the following way \cite{Bishop.2006}:
\begin{equation}\label{eq:errorfunctionNN}
\begin{split}
L_{error}(Y,\hat{Y})=\frac{1}{N}\sum_{k=1}^{N}\left|\left|Y_k-\hat{Y}_k\right|\right|_2^2.
\end{split}
\end{equation}
Depending on the shape of data, other norms than the mean squared error may be applied equally, e.g. the mean absolute error. However, the choice of the norm for the loss function \eqref{eq:errorfunctionNN} as well as the choice of activation functions and the number of neurons and layers strongly affect the network's performance and therefore should be chosen carefully and problem-specifically. Due to these numerous hyperparameters, influencing directly the training success, it is necessary to find the best suitable set of hyperparameters. Among different methods, e.g. random search or grid search, Bayesian optimisation has become the state-of-the-art technique to determine optimal hyperparameters for NNs. Eventually the error, resulting from the comparison \eqref{eq:errorfunctionNN}, is processed via backpropagation and used to adapt the neurons for a better fit \cite{Bishop.2006}. The capability of a NN as a global approximator of any function results from its structure, as it can be shown that it is confident to approximate any nonlinear mapping provided it holds a sufficient number of layers and neurons  \cite{Chen.1995}. Despite their success in many fields NNs usually lose physical interpretability, since they map inputs to outputs by nonlinear matrix operations, which cannot be interpreted for nonlinear dynamics. Moreover, the training of a NN is time-consuming and not straight-forward, since they tend towards overfitting \cite{Bishop.2006}. Therefore, data for training needs to be chosen carefully and should ideally cover most of the phase and input space to let the NN learn all possible situations and excitations, since it poorly extrapolates beyond the scope of training data. This leads to the need of numerous data samples for learning the dynamics of a system, which often is not possible to accomplish.

\subsection{Hybrid Modelling}\label{subsec:hybrid}
\begin{figure*}[!h]
	\centering
	\subfloat[Training of a PGNN(-L)]{\includegraphics[width=0.49\textwidth,scale=0.7]{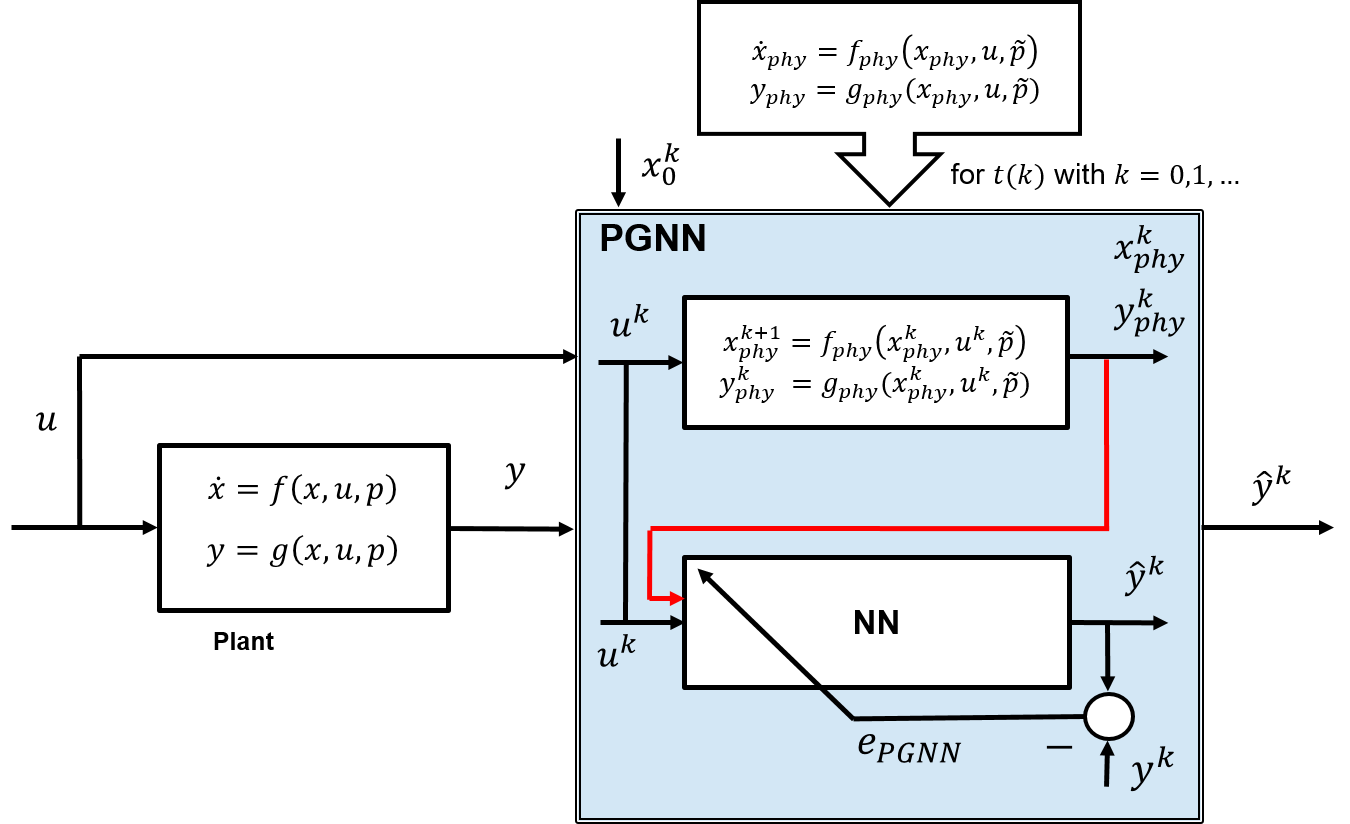}\label{fig:PGNNandPlant}}
	\hfill
	\subfloat[Prediction with a PGNN(-L)]{\includegraphics[width=0.49\textwidth,scale=0.7]{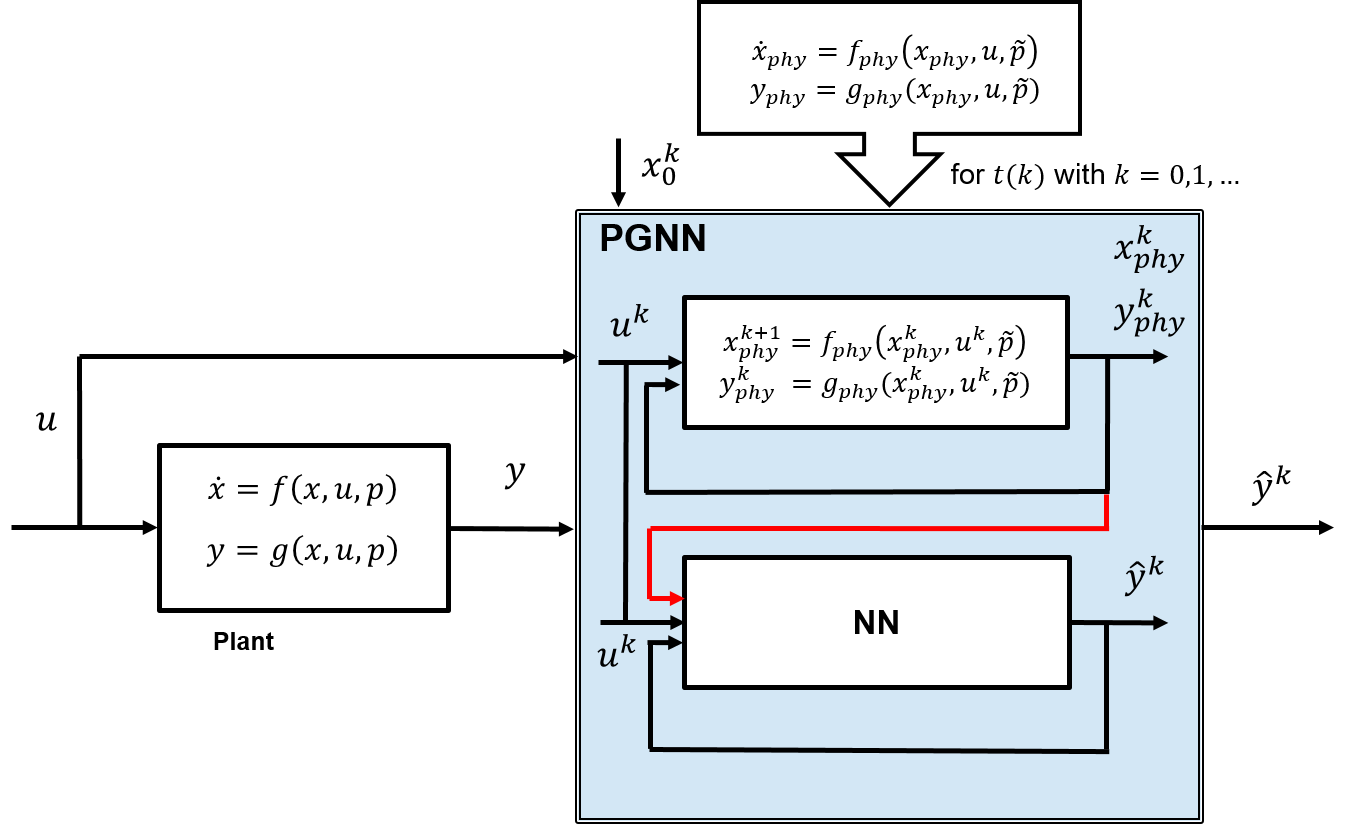}\label{fig:PGNNandPlantb}}
	\caption{Workflow of a plant and a PGNN(-L): (\ref{fig:PGNNandPlant}) During training the PGNN(-L) takes the excitation $u^k$, the simulated output $(x_{phy}^k,y_{phy}^k)$ of the discrete physical model and the measurement data $y^k$ into account to minimise the error $e_{PGNN}$ of its loss function \eqref{eq:errorfunctionNN} (or in case of PGNN-L \eqref{eq:physLoss}). (\ref{fig:PGNNandPlantb}) Once the PGNN(-L) is trained, the prediction of the previous time step $\hat{y}^k$ is constantly fed back into its input layer. With the current excitation $u^k$ and simulated output $(x_{phy}^k,y_{phy}^k)$ a prediction for the next time step can be estimated.}
\end{figure*}

After discussing techniques of purely physical and data-driven modelling in the previous sections, we now derive methods with both components, using prior physical knowledge and the availability of data to obtain a more accurate model of the considered system \eqref{eq:system}.\newline
As mentioned in the introduction, some hybrid approaches can be assigned to the group of \textit{Extended Residual Modelling}, which share the procedure of relating the model error to a prior physical model instead of only approximating the model error itself (Residual Modelling). \textit{Physics-Guided Neural Networks (PGNNs)} \cite{Karpatne.2017}, also known as \textit{Physics-Informed Neural Networks (PINNs)} \cite{Raissi.2021}, can be accounted to these techniques. They combine physics-based models and measurement data by either augmenting the neural network's inputs with the simulated outputs of the physics-based model or by utilising ODEs, which describe the system's dynamics, within the error function in the learning process.\newline
In this paper we adopt these techniques for non-autonomous systems, which we consider in control engineering,  and propose a method based on PGNNs and physical loss functions that copes with inconsistent physical models and noisy measurement data, whilst ensuring physical meaningful predictions. In figure \ref{fig:PGNNandPlant} the structure of a PGNN in relation with a plant during training is shown. The PGNN consists of a physical model \eqref{eq:physModel}, which is assumed to be not completely reliable, and is transformed into a discrete model, which can be evaluated by common integration techniques, such as explicit Euler or Runge-Kutta scheme. The PGNN's second component is a NN whose inner structure does not differ from a standard NN and can be chosen in a suitable manner. This includes e.g. feed-forward or recurrent connections of the neurons, the number of layers and neurons or the manner of activation functions, which all individually affect a network's performance. Both PGNN's components receive the current input signal $u^k$ the plant is stimulated with. The NN takes additionally the simulated output $(x_{phy}^k,y_{phy}^k)$ of the physical model \eqref{eq:physModel} and the measurement data $y^k$ into account to minimise the error $e_{PGNN}$ of its loss function and to provide a prediction $\hat{y}^k$. As NNs often tend to overfitting, an obligatory step of data preprocessing is the standardisation of data. Therefore, before using the excitation, simulated output or measurement data as inputs for the NN, this step has to be applied towards each of them. If data  feature different scales, normalisation or a logarithmic transformation of the samples may be considered additionally.\newline
In contrast to \cite{Cranmer.2020, Greydanus.2019, Yu.2020, Antonelo.2021} we take the prior physical knowledge as external input into account which rather guides the PGNN than forces it to strictly obey a priori presumed dynamics. Using such models of presumed dynamics within the loss function would imply we consider the model \eqref{eq:physModel} to be completely true. Due to the model-reality-gap this is hard to achieve and would lead the PGNN to learn possibly faulty dynamics. Subsequently, the PGNN would not yield a significant higher model accuracy than the presumed model and justify the effort of hybrid modelling or the loss of interpretability. Hence, it is more efficient to let the model \eqref{eq:physModel} guide the PGNN's learning process externally with a higher flexibility of adopting or neglecting information. 
However, data-driven methods in general and NNs in particular,  have significant drawbacks that result in poor solutions for unseen data or data beyond the training scope, a tendency for under- or overfitting and a lack of physically consistent predictions. The framework of PGNN-L overcomes these drawbacks and is completed with the incorporation of physically motivated loss functions during the learning procedure. Instead of only evaluating the traditional error loss function between the predicted and targeted output data, $\hat{Y}$ and $Y$ respectively, a second loss term is added with the weighting factor $\lambda_{phy}\in\left(0,1\right)$:
\begin{equation}\label{eq:physLoss}
\begin{split}
L(Y,\hat{Y},U)= (1-\lambda_{phy})L_{error}(Y,\hat{Y})+\lambda_{phy} L_{phy}(\hat{Y},U).
\end{split}
\end{equation} This may be an equality or inequality constraint $h:\,\mathbb{R}^l\mapsto \mathbb{R}$, taken into account e.g. by 
\begin{align}
\begin{split}
L_{phy}(\hat{Y},U)=\begin{cases}
\left|\left|h(\hat{Y},U)\right|\right|_2^2, &\text{for }h(\hat{Y},U)=0\\
\text{ReLu}(h(\hat{Y},U)), &\text{for }h(\hat{Y},U)\leq 0.
\end{cases}
\end{split}
\end{align}
Thus, physical or domain-specific constraints, derived from some understanding or experience with the plant, can be enabled. These could include conservation laws \cite{Raymond.2021}, like momentum or energy, any other physical relations \cite{Karpatne.2017}, laws of nature, feedback control laws, bounds \cite{Muralidhar.2018} or even vague knowledge, such as experts' experience or rudimentary rules. Extending the traditional loss function with a physically meaningful relation, the training itself is more efficient since the solution space is reduced. Furthermore, which is the major benefit, it enforces the learning process to strive for physically consistent solutions, which is considered as a robustness and safety criterion towards numerical instabilities or external disturbances if the constraint is chosen well. In this contribution we insert and investigate such a physical loss function for non-autonomous systems, which has not been covered yet. Since the control input $u$ means another degree of freedom it is even more difficult to learn and enforce physical consistency within a PGNN. We show that including the control input in the loss function \eqref{eq:physLoss} and respecting the energy conservation principle helps to guarantee physical consistency to a higher extent than other methods are able to guarantee. Another benefit of the physics-based loss function is its independence on the targeted outputs. This allows to continue the utilisation of it even after training, when no labelled data are available any more, and therefore ensures the constant check for physical consistency.  Nevertheless, the optimisation problem of minimising \eqref{eq:physLoss} is thus transformed into a multiobjective problem with competing loss criterions, which introduces new challenges, e.g. choosing a suitable method for finding the pareto front or pareto optimal solutions. Just recently, this challenge has been investigated in \cite{Rohrhofer.2021} for the first time. The authors demonstrate the effects of arbitrarily weighting different loss terms on the convergence of the learning process as well as on the possibly poor predictions. They conclude that strategies for weighting these loss terms need to be further examined. Nonetheless, for the scope of this work, we ensure that \eqref{eq:physLoss} yields a global minimum by using a simple method of scaling techniques, the weighted sum. Both the chosen error loss term and physical loss term are formulated by the mean squared error for the experiments in section \ref{sec:experiments}, resulting in a convex total loss function. Yet, we also examine the effects of the weighting factor to obtain insights into this novel-arised problem when balancing physics-based and traditional error loss terms.\newline
Using Algorithmic Differentiation \cite{Griewank.2008} the PGNN-L with multiple loss terms can be trained via backpropagation by the ADAM solver \cite{Kingma.2015}. Similar to the hyperparameter optimisation of standard NNs, the weighting factor $\lambda_{phy}$ is determined by Bayesian optimisation \cite{Shahriari.2016}, whose effects on the training and the prediction quality is explored in detail in section \ref{sec:experiments}. Once the PGNN-L is trained completely, it can be utilised for prediction as it is illustrated in figure \ref{fig:PGNNandPlantb}. It shows the prediction stage as an integration process which needs the current excitation $u^k$ as only information from the plant, the prediction $\hat{y}^k$ of the previous time step, which is constantly fed back into the NN as initial value, and the estimated prediction of the discrete physical model of the current time step, which is equally fed back into the physical model as initial value for the next time step. Provided the predictions of the simulation model \eqref{eq:physModel} are already calculated and the NN features a low complexity, it is very efficient and fast to evaluate the PGNN-L due to the computation of nonlinear matrix-vector-computations. However, even with simulating the physical model beforehand, it may be capable for real-time applications, depending on the number of neurons and layers as well as on the complexity of the physical model.  

\section{Experimental Results}\label{sec:experiments}
The previous described method is eventually demonstrated on two non-autonomous, nonlinear systems. From now on we assume that the collected data samples include the full state of the considered system ($y=x$), which is reasonable as an initial observer can be built upon prior knowledge of the physical model \eqref{eq:physModel}. This might lead to slightly corrupted data but can be compensated by the hybrid approach. For each considered system the following modelling approaches are compared throughout the paper: The typical physics-based model description, a solely data-driven, feed-forward neural network (NN), a PGNN with the same inner structure like the NN and a physical loss term (PGNN-L) as well as a SINDYc model based on the prior knowledge the physical model contains.
\subsection{Golf robot}  
\begin{figure}[h!]
	\subfloat[Construction]{\includegraphics[width=0.125\textwidth]{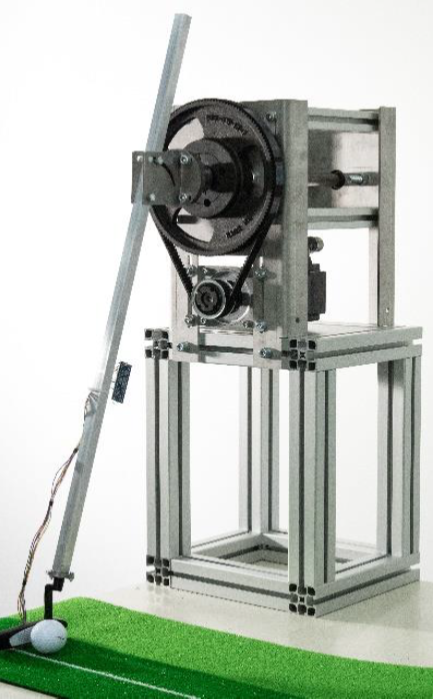}\label{fig:golfRobot}} \subfloat[Parameters $\tilde{p}$]{\begin{tabular}[b]{lrl}
			\hline
			\textbf{Parameter} & \textbf{Value} &\\
			\hline
			mass of club head $m$ & 0.5241&kg\\
			length from center of &&\\
			gravity to rotation axis $a$ &0.4702&m\\
			moment of inertia $J$ & 0.1445&kgm$^2$\\
			damping constant $d$ & 0.0132&kgs$^{-1}$\\
			distance friction point &&\\
			to rotation axis $r$ & 0.0245&m\\
			friction coefficient $\mu$ & 1.5136&-\\
			gravity $g$ & 9.81&ms$^{-2}$\\
			\hline
		\end{tabular}\label{tab:golfParameters}}
	\caption{Golf robot}
\end{figure}

\begin{figure*}[h!]
	\centering
	\includegraphics[width=\textwidth]{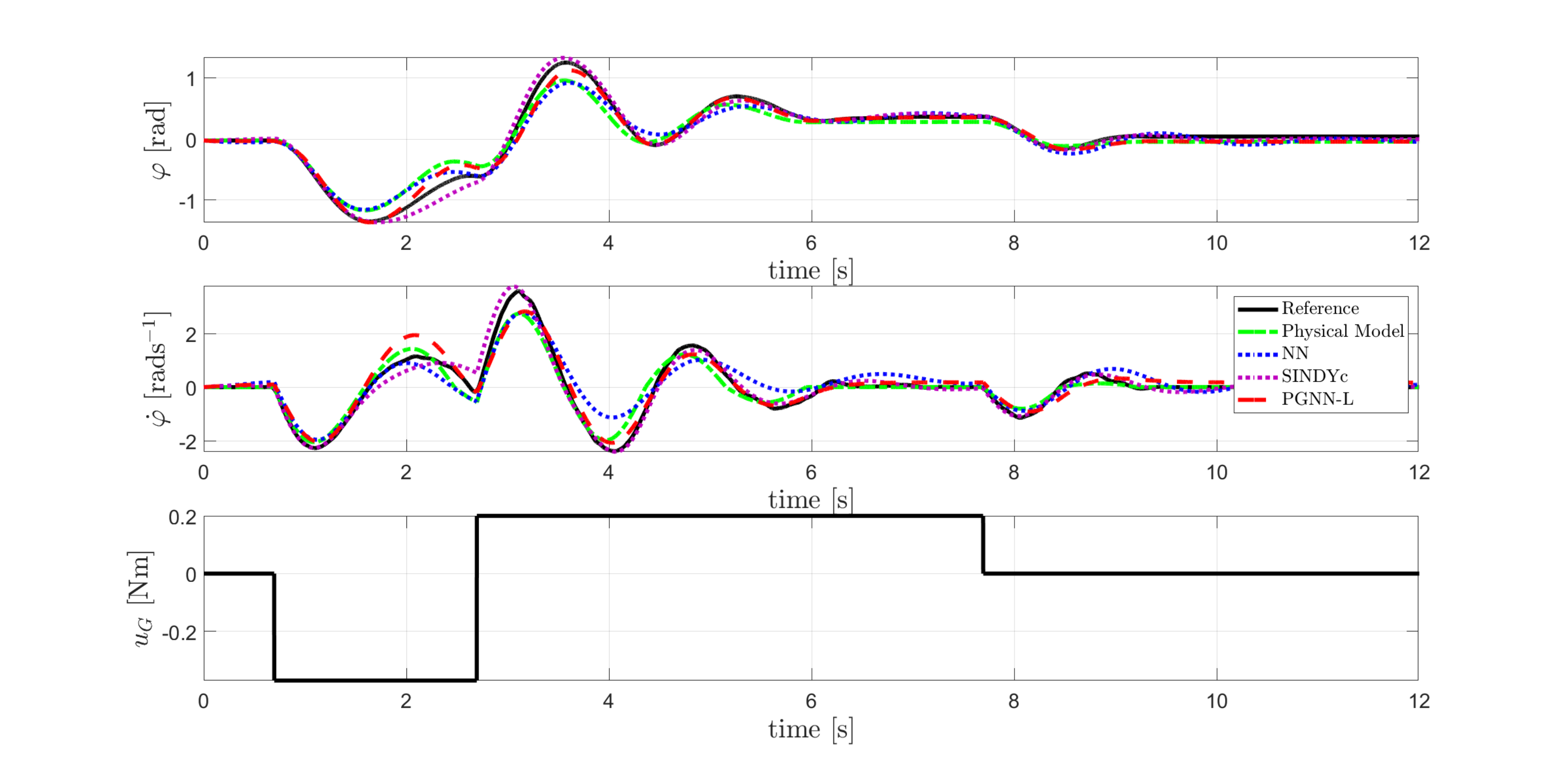}
	\caption{Prediction of golf robot's dynamics with $\lambda_{phy}=0.8175$, number of PGNN's and NN's neurons: 27 vs. 39}\label{fig:golfplot1}
\end{figure*}

The first system to be considered is a golf robot, a test rig at our laboratory to investigate machine learning based methods for hitting a golf ball within the green, which is illustrated in figure \ref{fig:golfRobot}. The physical model is captured by a nonlinear, input affine system of ODEs:
\begin{align}\label{eq:golfSystem}
\begin{split}
\dot{x} &=\begin{pmatrix}
x_2\\ \frac{-mga\sin x_1- F_G(x)+4u_G}{J}
\end{pmatrix},\\
F_G(x) &= d x_2+r\mu\, \text{sign}(x_2)\left(
mx_2^2 a+mg\cos x_1\right).
\end{split}
\end{align}
The motor torque, which is multiplied by four due to the gear ratio, serves as control input $u_G\in\mathbb{R}$, while the angle $\varphi$ and angular velocity $\dot{\varphi}$ of the club head describe the golf robot's states $x\in\mathbb{R}^2$. As stated before, only $\varphi$ is measured, whereas $\dot{\varphi}$ is estimated. The term $F_G(x)$ contains the dynamics resulting from static and sliding friction. The parameters $\tilde{p}=\left(m, a, J, d, r, \mu\right)^T \in\mathbb{R}^6$, estimated by measuring and a particle swarm optimisation, are shown in table \ref{tab:golfParameters}. Due to the stick-slip-effect the simulation of the model \eqref{eq:golfSystem} exhibits significant deviations compared to measurements from the test rig, which is illustrated in figure \ref{fig:golfplot1} by the green dashed and black solid lines, respectively.
To overcome this discrepancy, a PGNN-L with a simple feed-forward architecture, including two hidden layers, is trained by minimising the following specific form of the loss function \eqref{eq:physLoss}:
\begin{align}\label{eq:golfLoss}
\begin{split}
L(\hat{X},X,U)=&\left(1-\lambda_{phy}\right)\frac{1}{N}\sum_{k=1}^{N}(\hat{X}_k-X_k)^2\\
+\lambda_{phy}&\frac{1}{N-1}\sum_{k=2}^{N} \Delta E(\hat{X}_k,\hat{X}_{k-1},U_{k-1})^2\\
\Delta E(\hat{X}_k,\hat{X}_{k-1},U_{k-1})=&\Delta E_{kin}(\hat{X}_k,\hat{X}_{k-1})\\&+\Delta E_{pot}(\hat{X}_k,\hat{X}_{k-1})\\ &+\Delta E_{con}(\hat{X}_k,\hat{X}_{k-1},U_{k-1})\\&-\Delta E_{diss}(\hat{X}_k,\hat{X}_{k-1}),
\end{split}
\end{align}  
\begin{figure*}[h!]
\centering
\includegraphics[width=\textwidth]{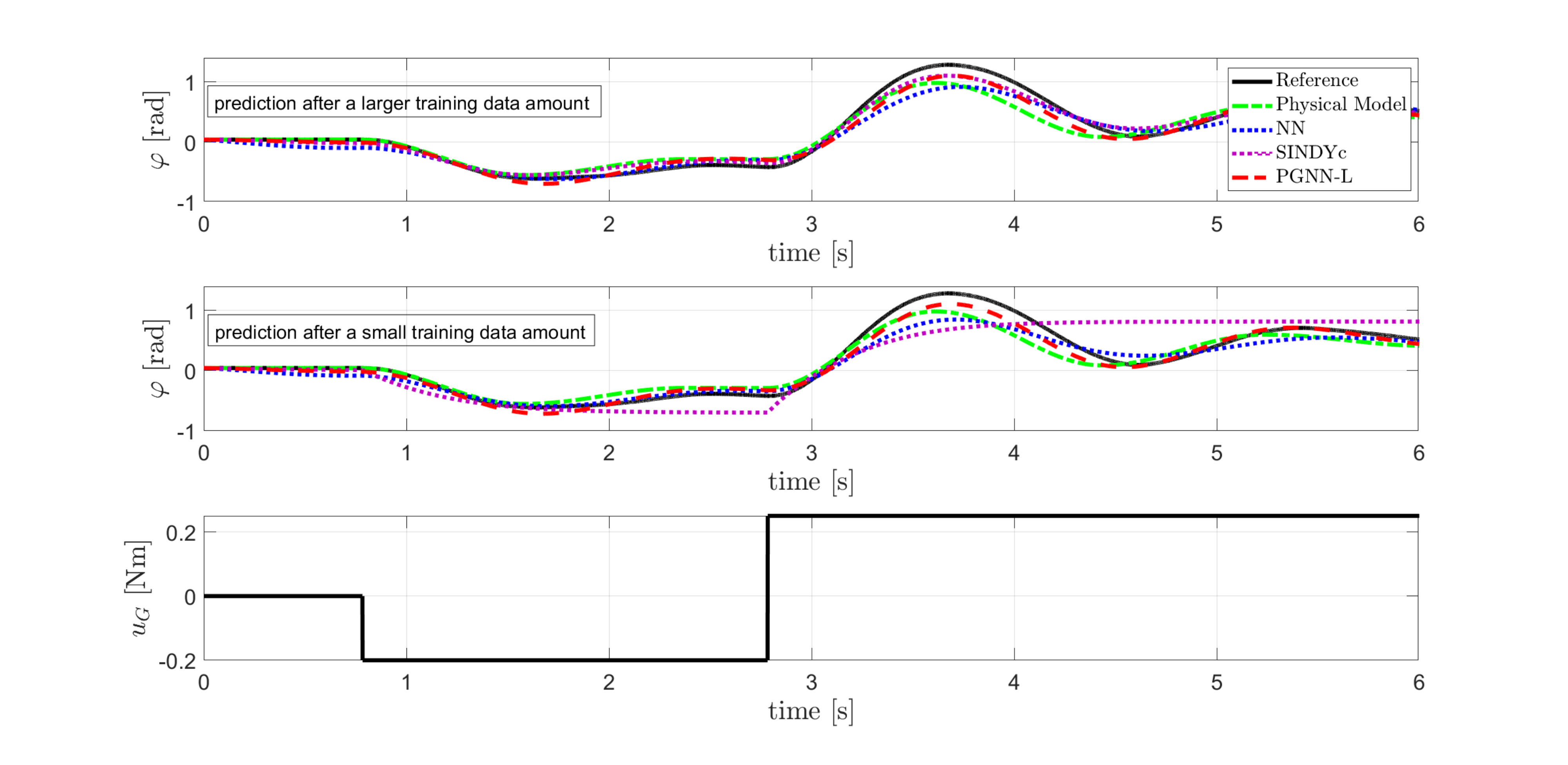}
\caption{Effects of the training data used: If only a small data amount covering the transient process is used, the SINDYc model will fail to reproduce the golf robot's angle due to missing information within the data. In contrast, the PGNN-L will capture the dynamics because of its additional information from the prior model and physical loss term.}\label{fig:datagolf}
\end{figure*}
for which the indices under each $\Delta E$ mark the kinetic, potential, control supplied and dissipative energy fractions, respectively. In contrast to \cite{Raymond.2021}, the second loss term extends the scope of conservative systems towards a variety of systems by taking additionally the dissipative and control supplied energy fractions into account. Hence, the second term in the loss function \eqref{eq:golfLoss} measures the energy balance between the previous and the current time step, which would be zero if the system's energy balance is modelled correctly. However, recalling that the prior knowledge and thus the energy balance are incomplete due to unmodelled partial dynamics, the physical loss term is independent of labelled data and is only evaluated on predicted outputs. In \cite{Raymond.2021} the authors compare the conservative system's energy balance between targeted and predicted outputs to ensure physical consistency. This is reasonable if the system is completely known. Yet, in this work this is not the case so the loss term evaluated on predicted outputs ensures the learning problem shifts to the direction of physically consistent solutions.\newline 
Data for training include measurement samples from six si\-nusoi\-dal-, step- and chirp-excited trajectories with a sampling rate of 1 kHz. These are separated in 60-20-20 percentage fractions for training, validation and test data sets. For both, the NN and the PGNN-L, a hyperparameter optimisation is performed by Bayesian optimisation to obtain the best set of affecting parameters, including $\lambda_{phy}$ for the PGNN-L. The SINDYc library vector $\Psi$ contains the states, the control input as well as trigonometric functions to find a suitable $\Xi$. Eventually, a reference trajectory from the test rig is used for simulation to quantify the models' accuracy, which is evaluated between the measured angle $\varphi^k$ and the corresponding prediction $\hat{\varphi}^k$.\newline
Figure \ref{fig:golfplot1} shows the resulting trajectories of each considered method. It reveals that the PGNN-L is confident to capture the golf robot's dynamics to a higher accuracy than the prior model. This is verified by the error displayed in table \ref{tab:error}. Although it covers the dynamics already very closely and does not exhibit significant phase shifts, it sometimes fails to reproduce the exact amplitude. However, it outperforms a traditional NN clearly, which is even worse than the physical model. For instance, the NN is not capable to recognise the friction and still oscillates a long time after the system has reached its equilibrium. Furthermore, the NN's complexity is much higher than the PGNN-L's, which only needs 27 neurons. Compared to the SINDYc model, the PGNN-L's accuracy appears to be slightly worse which is confirmed by the error in table \ref{tab:error}. Nevertheless, if training data only covers a small region of the phase and input space, e.g. the transient process of each of the previous for training used six trajectories, the PGNN-L approach clearly is superior over SINDYc. This is illustrated in figure \ref{fig:datagolf}, whose top figure shows simulated trajectories of the angle after training with data that covers a reasonable amount of the phase and input space and whose medium figure displays results of the angle after training with data that covers only the transient process. In this situation the advantage of the PGNN-L becomes clear: Since it relies on the physical model and a qualitative loss term, it can cope with incomplete data in contrast to SINDYc which approximates the dynamics poorly. Due to its reliance on data and a few governing terms SINDYc needs representative data of the entire phase and input space to approximate a confident model. Therefore, if only incomplete data combined with some physical understanding is available, which is often the case in complex industrial applications, PGNN-L is the optimal choice to derive a model.\newline 
\begin{figure}[b]
	\centering
	\includegraphics{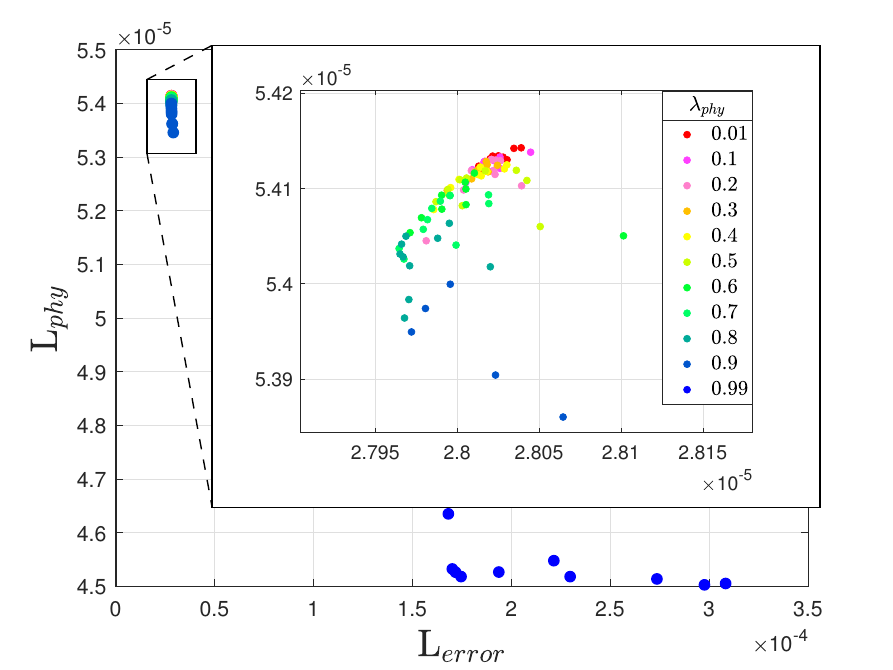}
	\caption{Pareto front for varying weighting factors $\lambda_{phy}$ in equation \eqref{eq:physLoss}}\label{fig:paretofront}
\end{figure}
Besides the PGNN-L's accuracy the trade-off between the competing loss terms $L_{error}$ and $L_{phy}$ during training is examined. In figure \ref{fig:paretofront} the pareto front is shown for varying weighting factors $\lambda_{phy}$. If $\lambda_{phy}$ is 0.99, which means the focus on the data-fit term \eqref{eq:errorfunctionNN} is very  low, the error between targeted and predicted angle is high. However, zooming into the rectangle reveals a typical convex form of a pareto front. It indicates that a good weighting factor for the golf robot is between 0.8 and 0.9 meaning the two competing loss terms are both very small. This is a high emphasis on the energy conservation loss term. Although the PGNN-L model is at least as good as the SINDYc model or in case of low data availability even more precise, simulations as in figure \ref{fig:golfplot1} and figure \ref{fig:datagolf} still exhibit some deviations between the PGNN-L model and the true dynamics, e.g. in the amplitudes. Therefore, a better fit towards the data would be desirable. Instead of setting $\lambda_{phy}$ constant overall the learning process, $\lambda_{phy}$ may be regarded as an adaptive parameter similar to the learning rate which is able to vary according to the training loss. Future work will examine this strategy.

\subsection{Servo valve}
\begin{figure}[h]
	\centering
	\includegraphics[scale=0.15]{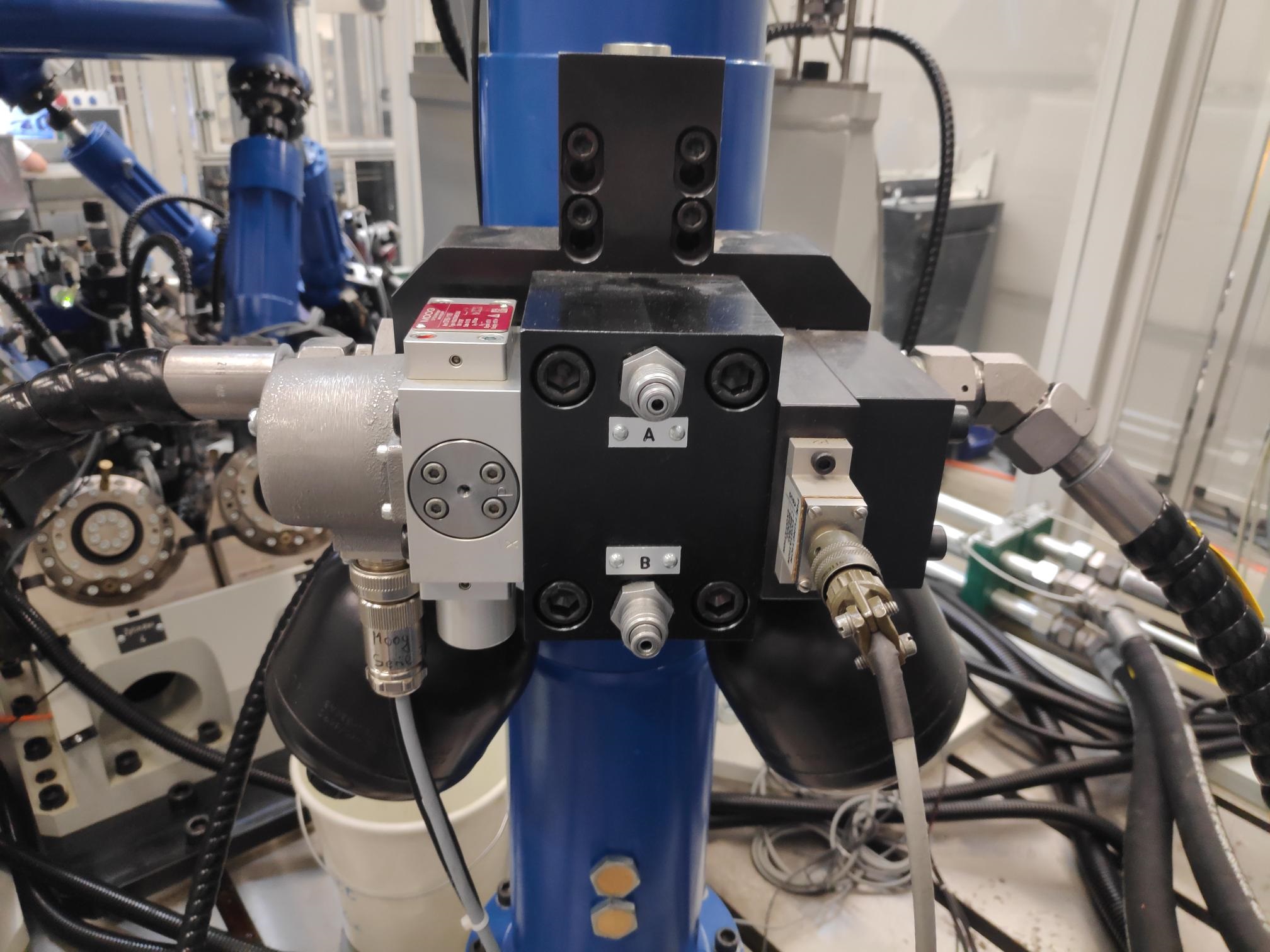}
	\caption{Servo valve deployed in a test rig}\label{fig:valve}
\end{figure}

In many hydraulic applications servo valves are frequently utilised as actuators. However, servo valves often feature complex nonlinear dynamics, for instance the bandwidth of the valve strongly depends on the amplitude of the input signal. Therefore, it is common to utilise simplified models for valve dynamics. An approach with very low modelling depth is e.g. the description of the servo valve dynamics via a second order system:
\begin{equation}\label{eq:valveSystem}
\dot{x}=\begin{pmatrix}
x_2\\ -\frac{2D_V}{T_V} x_2 -\frac{1}{T_V^2} x_1 + \frac{K_V}{T_V^2}u_V
\end{pmatrix},
\end{equation}
with $x\in\mathbb{R}^2$ describing the position $y_V$ 
and the velocity $\dot{y}_V$ of the valve slider, $u_V\in\mathbb{R}$ representing the voltage the servo valve receives as input signal and $T_V,D_V$ and $K_V$ denoting the typical parameters of a second order system. Again, only $y_V$ is measured, whereas $\dot{y}_V$ is estimated. In this work we consider a four way three position servo valve consisting of a valve slider, a torque motor and a hydraulic amplifier based on the nozzle-baffle principle which is illustrated in figure \ref{fig:valve}.  
\begin{table}[h]
	\caption{Parameters $\tilde{p}$}\label{tab:valveParameters}
	\centering
	\begin{tabular}[b]{lrl}
		\hline
		Parameter & Value &\\
		\hline
		natural frequency $f_V$ & 350 &Hz\\
		damping ratio $D_V$ & 0.5 &-\\
		gain factor $K_V$ & 0.1 &- \\
		voltage $u_V$ & -10 $\dots$ 10 &V\\
		maximal position $y_{V,max}$ & 4.2672$\cdot$10$^{-4}$ &m\\
		supply pressure $p_S$ & 280 &bar\\
		\hline
	\end{tabular}
\end{table}
\begin{figure*}[h!]\label{fig:valveplot}
	\centering
	\subfloat[PGNN-L with prior model that does not respect physical limitations]{\includegraphics[width=0.5\textwidth]{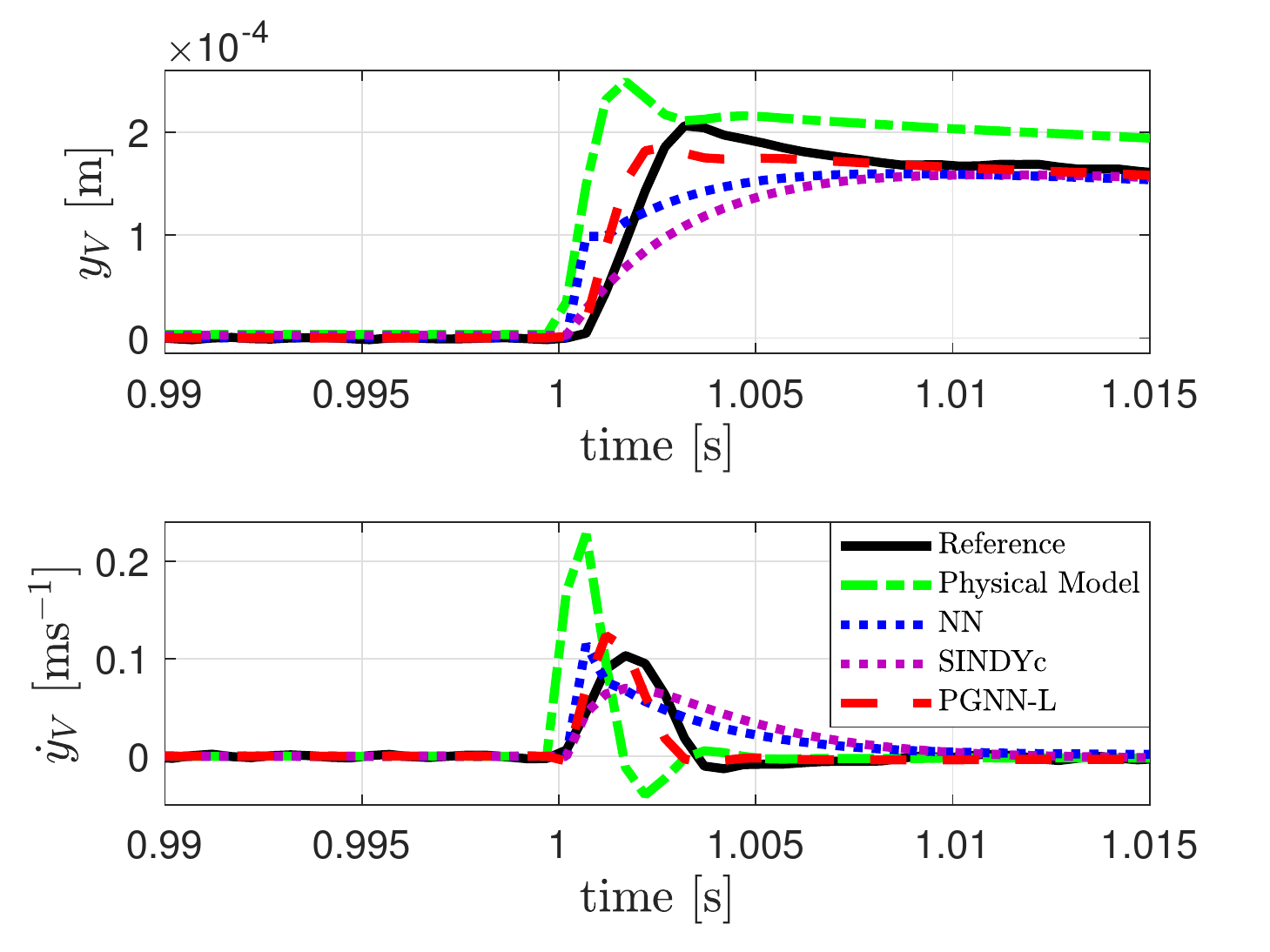}\label{fig:valveplot1}}
	\subfloat[PGNN-L with prior model that respects physical limitations]{\includegraphics[width=0.5\textwidth]{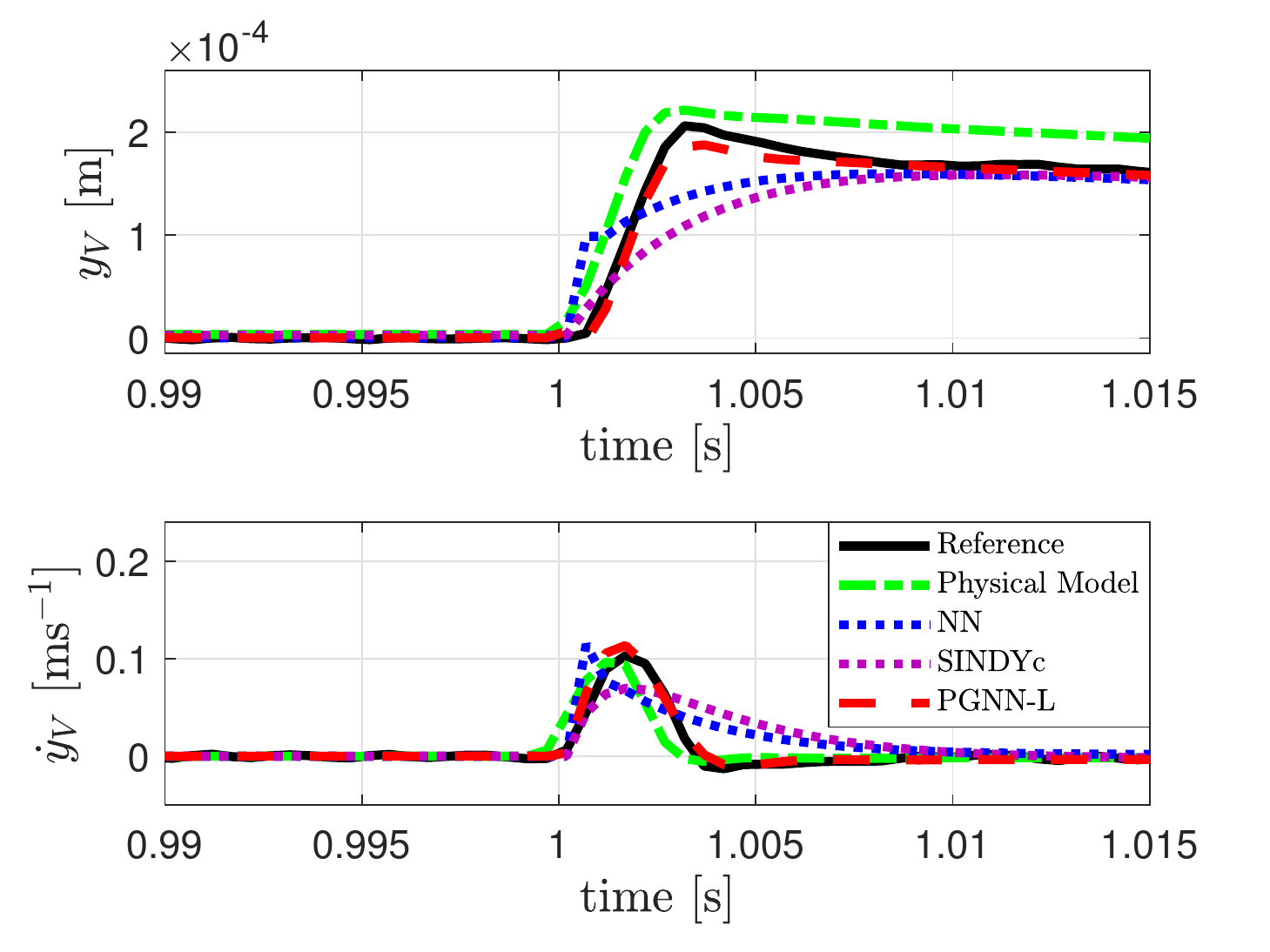}\label{fig:valveplot2}}
	\caption{Prediction of servo valve dynamics by a step excitation of 5$V$ at 1$s$ : (\ref{fig:valveplot1}) $\lambda_{phy}=0.2527$, number of PGNN-L's and NN's neurons: 11 vs. 92, (\ref{fig:valveplot2})  $\lambda_{phy}=0.3206$, number of PGNN-L's and NN's neurons: 11 vs. 92}
\end{figure*}

The valve slider is controlled with electric position control. The servo valve's parameters for the physical model \eqref{eq:valveSystem} are displayed in table \ref{tab:valveParameters}.
Although the model \eqref{eq:valveSystem} is able to predict the system's behaviour in general, it fails to preserve the intrinsic limitations of the valve slider's velocity and acceleration due to the oil viscosity inside the piston and due to the torque motor's position because of its operational limit. This is illustrated in figure \ref{fig:valveplot1} by the green dotted and the black solid lines, respectively. Analogue to the first considered system, a PGNN-L with feed-forward structure and two hidden layers is trained via the loss function \eqref{eq:golfLoss}.
Data for training include ten, slightly noisy trajectories of different step excitations $u_V$ with a sampling rate of $2$ kHz and are separated in the same manner (60-20-20) for training, validation and test data sets as in the previous section. Furthermore, for both NN and PGNN-L the best set of hyperparameters is determined by Bayesian optimisation, including $\lambda_{phy}$ for the PGNN-L. The library vector for SINDYc contains the states $y_V$, $\dot{y}_V$ and the voltage $u_V$ within this example. Each model's accuracy is quantified by the same reference trajectory from measurement data, which has not been covered during training, and evaluated between the measured state $y_V^k$ and the predicted state $\hat{y}_V^k$.\newline
Figure \ref{fig:valveplot1} shows that the PGNN-L is able to predict the servo valve's behaviour to a higher accuracy than the compared methods are able to achieve. Although the physical model reveals significant deviations especially in the magnitude, e.g. in the slider's velocity, the PGNN-L can process this information and utilize it with data and the energy respecting loss term to capture the true servo valve's dynamics closely. Neither the standard NN nor the SINDYc model can reproduce the dynamics to this accuracy. Both fail in recognising the oscillating and damping effects. If the physical model \eqref{eq:valveSystem} is improved by estimated bounds of the slider's velocity and acceleration due to the oil's viscosity, this information can be preserved by the PGNN-L according to figure \ref{fig:valveplot2}. Here, the PGNN-L reproduces the valve's real-word dynamics even more precisely compared to the setting before. A much lower error in table \ref{tab:error} verifies this observation. Hence, the modelling depth of the prior model appears to directly affect the PGNN-L's performance and indicates a specific modelling depth is essential for deriving meaningful hybrid models. However, the number of neurons remains constant which need not necessarily mean the modelling depth does not influence the PGNN's model complexity. One would expect the more precise the prior model is the better the PGNN-L's approximation of the true dynamics will be. This needs to be examined further in future work. Nonetheless, the PGNN-L is less complex than a standard NN and features only an eighth of the NN's number of neurons. This benefit is caused by the utilisation of the prior model and the physically motivated loss term. Moreover, it is remarkable that the SINDYc model is not capable to recognise any intrinsic bounds and therefore reproduces the valve's dynamics very poorly. Since SINDYc operates by providing governing terms, it cannot respect any limitations. Thus, the PGNN-L is the optimal choice to cover the valve's real-world dynamics among the examined techniques.       
\begin{table}[h]	
	\caption{Root-mean-squared error (RMSE) of the considered nonlinear systems}
	\begin{tabular}{lccc}
		\hline
		Method & Golf robot (fig. \ref{fig:golfplot1}) & Valve (fig. \ref{fig:valveplot1}) & Valve  (fig. \ref{fig:valveplot2})\\
		\hline
		PhysModel & 6.6563$\cdot$ 10$^{-4}$ & 5.0501$\cdot$ 10$^{-4}$ & 3.9959$\cdot$ 10$^{-4}$\\
		NN & 8.9736$\cdot$ 10$^{-4}$ & 2.1382$\cdot$ 10$^{-4}$ & 2.1382$\cdot$ 10$^{-4}$\\
		SINDYc &  3.4516$\cdot$ 10$^{-4}$ & 2.8619$\cdot$ 10$^{-4}$& 2.8619$\cdot$ 10$^{-4}$\\
		PGNN-L &  3.7982$\cdot$ 10$^{-4}$ & 1.1982$\cdot$ 10$^{-4}$ & 7.3835$\cdot$ 10$^{-5}$\\
		\hline
	\end{tabular}
	\label{tab:error}\end{table}

\section{Conclusion and outlook}\label{sec:conclusion}
In this contribution we extend the framework of PGNNs towards system identification of non-autonomous systems in control engineering. The described approach PGNN-L consists of a PGNN with augmented inputs by the control input and simulated output of the corresponding prior model as well as of a physics-motivated loss term within the learning procedure. Two real-world nonlinear systems are examined, for which the physical loss term respects the energy conservation principle and takes also dissipative and control supplied energy fractions into account. It is shown that the method achieves a higher model accuracy than existing methods, such as NN or SINDYc models, if incomplete or noisy data is available. It also respects physical limitations. Additionally, it enables mostly lower model complexity, since it is guided by the prior physical model. Furthermore, its physics-motivated loss term ensures physically consistent predictions to a higher extent, which is considered as robustness criterion.
Future work include the investigation of balancing more competing physical loss terms, including vague knowledge, the effects of different modelling depths for the prior model and the development of an adaptive weighting factor schedule for the physics-based loss term.

\section*{Acknowledgment}

This work was developed in the junior research group DART (Daten\-ge\-trie\-be\-ne Methoden in der Regelungstechnik), University Paderborn, and funded by the Federal Ministry of Education and Research of Germany (BMBF - Bundesministerium für Bildung und Forschung) under the funding code 01IS20052. The responsibility for the content of this publication lies with the authors.


\begin{thebibliography}{10}
	\providecommand{\url}[1]{#1}
	\csname url@samestyle\endcsname
	\providecommand{\newblock}{\relax}
	\providecommand{\bibinfo}[2]{#2}
	\providecommand{\BIBentrySTDinterwordspacing}{\spaceskip=0pt\relax}
	\providecommand{\BIBentryALTinterwordstretchfactor}{4}
	\providecommand{\BIBentryALTinterwordspacing}{\spaceskip=\fontdimen2\font plus
		\BIBentryALTinterwordstretchfactor\fontdimen3\font minus
		\fontdimen4\font\relax}
	\providecommand{\BIBforeignlanguage}[2]{{%
			\expandafter\ifx\csname l@#1\endcsname\relax
			\typeout{** WARNING: IEEEtran.bst: No hyphenation pattern has been}%
			\typeout{** loaded for the language `#1'. Using the pattern for}%
			\typeout{** the default language instead.}%
			\else
			\language=\csname l@#1\endcsname
			\fi
			#2}}
	\providecommand{\BIBdecl}{\relax}
	\BIBdecl
	
	
	\bibitem{Brunton.2016}
	S.~L. Brunton, J.~L. Proctor, and J.~N. Kutz, ``Sparse identification of
	nonlinear dynamics with control (sindyc),'' in \emph{IFAC-PapersOnline},
	2016, vol.~49, pp. 710--715.
	
	\bibitem{Karpatne.2017}
	\BIBentryALTinterwordspacing
	A.~Karpatne, G.~Atluri, J.~Faghmous, M.~Steinbach, A.~Banerjee, A.~Ganguly,
	S.~Shekhar, N.~Samatova, and V.~Kumar, ``Theory-guided data science: A new
	paradigm for scientific discovery from data,'' in \emph{IEEE Transactions on
		Knowledge and Data Engineering}, 2017, vol.~29, pp. 2318--2331. [Online].
	Available: \url{http://arxiv.org/pdf/1612.08544v2}
	\BIBentrySTDinterwordspacing
	
	\bibitem{Psichogios.1995}
	D.~C. Psichogios and L.~H. Ungar, ``A hybrid neural network-first principles
	approach to process modeling,'' in \emph{AIChE Journal}, 1995, vol.~38.
	
	\bibitem{Willard.2020}
	\BIBentryALTinterwordspacing
	J.~Willard, X.~Jia, S.~Xu, M.~Steinbach, and V.~Kumar, ``Integrating
	physics-based modeling with machine learning: A survey,'' 2020. [Online].
	Available: \url{http://arxiv.org/pdf/2003.04919v4}
	\BIBentrySTDinterwordspacing
	
	\bibitem{Karpatne.2017b}
	\BIBentryALTinterwordspacing
	A.~Karpatne, W.~Watkins, J.~Read, and V.~Kumar, ``Physics-guided neural
	networks (pgnn): An application in lake temperature modeling,'' 2017.
	[Online]. Available: \url{http://arxiv.org/pdf/1710.11431v2}
	\BIBentrySTDinterwordspacing
	
	\bibitem{Antonelo.2021}
	\BIBentryALTinterwordspacing
	E.~A. Antonelo, E.~Camponogara, L.~O. Seman, E.~R. de~Souza, J.~P. Jordanou,
	and J.~F. Hubner, ``Physics-informed neural nets-based control,'' 2021.
	[Online]. Available: \url{http://arxiv.org/pdf/2104.02556v1}
	\BIBentrySTDinterwordspacing
	
	\bibitem{Pathak.2018}
	J.~Pathak, A.~Wikner, R.~Fussell, S.~Chandra, B.~R. Hunt, M.~Girvan, and
	E.~Ott, ``\BIBforeignlanguage{eng}{Hybrid forecasting of chaotic processes:
		Using machine learning in conjunction with a knowledge-based model},'' in
	\emph{\BIBforeignlanguage{eng}{Chaos}}, 2018, vol.~28, p. 041101.
	
	\bibitem{Yu.2020}
	Y.~Yu, H.~Yao, and Y.~Liu, ``Structural dynamics simulation using a novel
	physics-guided machine learning method,'' in \emph{Engineering Applications
		of Artificial Intelligence}, vol.~96, 2020, p. 103947.
	
	\bibitem{Yang.2020}
	\BIBentryALTinterwordspacing
	L.~Yang, X.~Meng, and G.~E. Karniadakis, ``B-pinns: Bayesian physics-informed
	neural networks for forward and inverse pde problems with noisy data,'' 2020.
	[Online]. Available: \url{http://arxiv.org/pdf/2003.06097v1}
	\BIBentrySTDinterwordspacing
	
	\bibitem{Muralidhar.2018}
	N.~Muralidhar, M.~R. Islam, M.~Marwah, A.~Karpatne, and N.~Ramakrishnan,
	``Incorporating prior domain knowledge into deep neural networks,'' in
	\emph{2018 IEEE International Conference on Big Data}, 2018, pp. 36--45.
	
	\bibitem{Cranmer.2020}
	\BIBentryALTinterwordspacing
	M.~Cranmer, S.~Greydanus, S.~Hoyer, P.~Battaglia, D.~Spergel, and S.~Ho,
	``Lagrangian neural networks,'' 2020. [Online]. Available:
	\url{http://arxiv.org/pdf/2003.04630v2}
	\BIBentrySTDinterwordspacing
	
	\bibitem{Greydanus.2019}
	\BIBentryALTinterwordspacing
	S.~Greydanus, M.~Dzamba, and J.~Yosinski, ``Hamiltonian neural networks,''
	2019. [Online]. Available: \url{http://arxiv.org/pdf/1906.01563v3}
	\BIBentrySTDinterwordspacing
	
	\bibitem{Raymond.2021}
	\BIBentryALTinterwordspacing
	S.~J. Raymond and D.~B. Camarillo, ``Applying physics-based loss functions to
	neural networks for improved generalizability in mechanics problems,'' in
	\emph{Journal of Computational Physics}, 2021. [Online]. Available:
	\url{http://arxiv.org/pdf/2105.00075v1}
	\BIBentrySTDinterwordspacing
	
	\bibitem{Rueden.2019}
	\BIBentryALTinterwordspacing
	L.~von Rueden, S.~Mayer, K.~Beckh, B.~Georgiev, S.~Giesselbach, R.~Heese,
	B.~Kirsch, J.~Pfrommer, A.~Pick, R.~Ramamurthy, M.~Walczak, J.~Garcke,
	C.~Bauckhage, and J.~Schuecker, ``Informed machine learning -- a taxonomy and
	survey of integrating knowledge into learning systems,'' 2019. [Online].
	Available: \url{http://arxiv.org/pdf/1903.12394v2}
	\BIBentrySTDinterwordspacing
	
	\bibitem{Lutter.2015}
	\BIBentryALTinterwordspacing
	M.~Lutter, C.~Ritter, and J.~Peters, ``Deep lagrangian networks: Using physics
	as model prior for deep learning,'' in \emph{ICLR}, 2015. [Online].
	Available: \url{http://arxiv.org/pdf/1907.04490v1}
	\BIBentrySTDinterwordspacing
	
	\bibitem{Daw.2019}
	\BIBentryALTinterwordspacing
	A.~Daw, R.~Q. Thomas, C.~C. Carey, J.~S. Read, A.~P. Appling, and A.~Karpatne,
	``Physics-guided architecture (pga) of neural networks for quantifying
	uncertainty in lake temperature modeling,'' 2019. [Online]. Available:
	\url{http://arxiv.org/pdf/1911.02682v1}
	\BIBentrySTDinterwordspacing
	
	\bibitem{Chen.2018}
	B.~Chen, H.~Zhang, X.~Liu, and C.~Lin, ``\BIBforeignlanguage{eng}{Neural
		observer and adaptive neural control design for a class of nonlinear
		systems},'' in \emph{\BIBforeignlanguage{eng}{IEEE Transactions on Neural
			Networks and Learning Systems}}, 2018, vol.~29, pp. 4261--4271.
		
	\bibitem{Li.}
	\BIBentryALTinterwordspacing
	Z.~Li, N.~Kovachki, K.~Azizzadenesheli, B.~Liu, K.~Bhattacharya, A.~Stuart, and
	A.~Anandkumar, ``Neural operator: Graph kernel network for partial
	differential equations,'' 2020. [Online]. Available:
	\url{http://arxiv.org/pdf/2003.03485v1}
	\BIBentrySTDinterwordspacing
	
	\bibitem{Li.2021}
	\BIBentryALTinterwordspacing
	------, ``Fourier neural operator for parametric partial differential
	equations,'' in \emph{ICLR}, 2021. [Online]. Available:
	\url{http://arxiv.org/pdf/2010.08895v2}
	\BIBentrySTDinterwordspacing
	
	\bibitem{Rohrhofer.2021}
	\BIBentryALTinterwordspacing
	F.~M. Rohrhofer, S.~Posch, and B.~C. Geiger, ``On the pareto front of
	physics-informed neural networks,'' 2021. [Online]. Available:
	\url{http://arxiv.org/pdf/2105.00862v1}
	\BIBentrySTDinterwordspacing
	
	\bibitem{Bishop.2006}
	C.~M. Bishop, \emph{Pattern recognition and machine learning}, ser. Information
	science and statistics.\hskip 1em plus 0.5em minus 0.4em\relax Springer,
	2006.
		
	\bibitem{Chang.}
	\BIBentryALTinterwordspacing
	Y.-C. Chang, N.~Roohi, and S.~Gao, ``Neural lyapunov control,'' 2019. [Online].
	Available: \url{https://arxiv.org/abs/2005.00611}
	\BIBentrySTDinterwordspacing
	
	\bibitem{Chen.2018b}
	\BIBentryALTinterwordspacing
	R.~T.~Q. Chen, Y.~Rubanova, J.~Bettencourt, and D.~Duvenaud, ``Neural ordinary
	differential equations,'' in \emph{32nd Conference on Neural Information
		Processing Systems (NeurIPS)}, 2018. [Online]. Available:
	\url{http://arxiv.org/pdf/1806.07366v5}
	\BIBentrySTDinterwordspacing
	
	\bibitem{Chen.1995}
	T.~Chen and H.~Chen, ``\BIBforeignlanguage{eng}{Approximation capability to
		functions of several variables, nonlinear functionals, and operators by
		radial basis function neural networks},'' in
	\emph{\BIBforeignlanguage{eng}{IEEE Transactions on Neural Networks}}, 1995,
	vol.~6, pp. 904--910.
	
	\bibitem{Raissi.2021}
	\BIBentryALTinterwordspacing
	M.~Raissi, P.~Perdikaris, and G.~E. Karniadakis, ``Physics-informed neural
	networks: A deep learning framework for solving forward and inverse problems
	involving nonlinear partial differential equations,'' in \emph{Journal of
		Computational Physics}, vol. 378, 2021, pp. 686--707. [Online]. Available:
	\url{https://www.sciencedirect.com/science/article/pii/S0021999118307125}
	\BIBentrySTDinterwordspacing
	
	\bibitem{Griewank.2008}
	A.~Griewank and A.~Walther, \emph{\BIBforeignlanguage{eng}{Evaluating
			derivatives}}, 2nd~ed.\hskip 1em plus 0.5em minus 0.4em\relax {SIAM Soc. f.
		Industrial and Applied Mathematics}, 2008, vol. 105.

	\bibitem{Kingma.2015}
	\BIBentryALTinterwordspacing
	D.~P. Kingma and J.~Ba, ``Adam: A method for stochastic optimization,'' in
	\emph{ICLR}, 2015. [Online]. Available:
	\url{http://arxiv.org/pdf/1412.6980v9}
	\BIBentrySTDinterwordspacing
	
	\bibitem{Shahriari.2016}
	B.~Shahriari, K.~Swersky, Z.~Wang, R.~P. Adams, and N.~de~Freitas, ``Taking the
	human out of the loop: A review of bayesian optimization,'' in
	\emph{Proceedings of the IEEE}, 2016, vol. 104, pp. 148--175.
		
\end{thebibliography}
\end{document}